\documentclass[12pt]{article}
\usepackage{amssymb}
\usepackage{dsfont,tikz,hyperref,url}
\newcommand{\N}{\mathds{N}}
\newcommand{\Z}{\mathds{Z}}

\newcommand{\R}{\mathds{R}}

\begin{document}

\begin{center}
{\bf Some 2-adic conjectures concerning \\
polyomino tilings of Aztec diamonds} \\
\ \\
James Propp, UMass Lowell \\
August 10, 2022 \\
\ \\
{\it Dedicated to Michael Larsen on the occasion of his 60th birthday}
\end{center}

\vspace{0.5in}

\begin{abstract}
\noindent
For various sets of tiles, we count the ways to tile 
an Aztec diamond of order $n$ using tiles from that set. 
The resulting function $f(n)$ often has interesting behavior 
when one looks at $n$ and $f(n)$ modulo powers of 2.
\end{abstract}

\begin{section}{Introduction}
\label{sec:intro}

I had a great time working on domino tilings of Aztec diamonds 
with Noam Elkies, Greg Kuperberg, and Michael Larsen back in the late 1980s, 
and the paper we wrote together~\cite{EKLP} 
had a huge impact on my career.
So I’d like to honor Michael by proposing 
some new problems about tilings of Aztec diamonds 
(and other regions), many of which are more challenging 
than the one I shared with him thirty-something years ago
and have a number-theoretic slant that I think he will enjoy.
Ideally the solutions to these problems will involve 
interesting applications of algebra to combinatorics.

\begin{figure}[h]
\begin{center}
\includegraphics[width=2in]{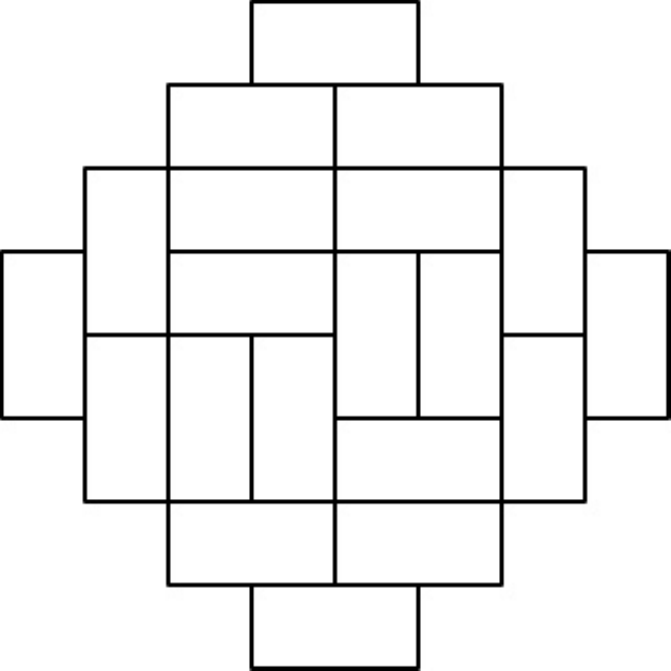}
\end{center}
\caption{A domino tiling of the Aztec diamond of order 4.}
\label{fig:dominos}
\end{figure}

Here is some general background.

The main result of~\cite{EKLP} was that the number of domino-tilings 
of an Aztec diamond of order $n$ is $2^{n(n+1)/2}$
(\href{https://oeis.org/A006125}{A006125}),
where a domino is a rectangle in $\R^2$ of the form
$[i,i+1] \times [j,j+2]$ or $[i,i+2] \times [j,j+1]$ (with $i,j \in \Z$)
and the Aztec diamond of order $n$  
is the union of the squares $[i,i+1] \times [j,j+1]$
lying entirely within the region $\{(x,y): \ |x|+|y| \leq n+1\}$.
Figure~\ref{fig:dominos} shows one of the $2^{(4)(5)/2}$ domino tilings 
of the Aztec diamond of order 4.

Mihai Ciucu~\cite{Ci} proved combinatorially
that the number of domino tilings of the $2n$-by-$2n$ square
(\href{https://oeis.org/A004003}{A004003})
can be written in the form $2^n f(n)^2$
where $f(n)$ is the number of domino tilings of the region
exemplified for $n=4$ in Figure~\ref{fig:ciucu}
(\href{https://oeis.org/A065072}{A065072}).

\begin{figure}[h]
\begin{center}
\includegraphics[width=2in]{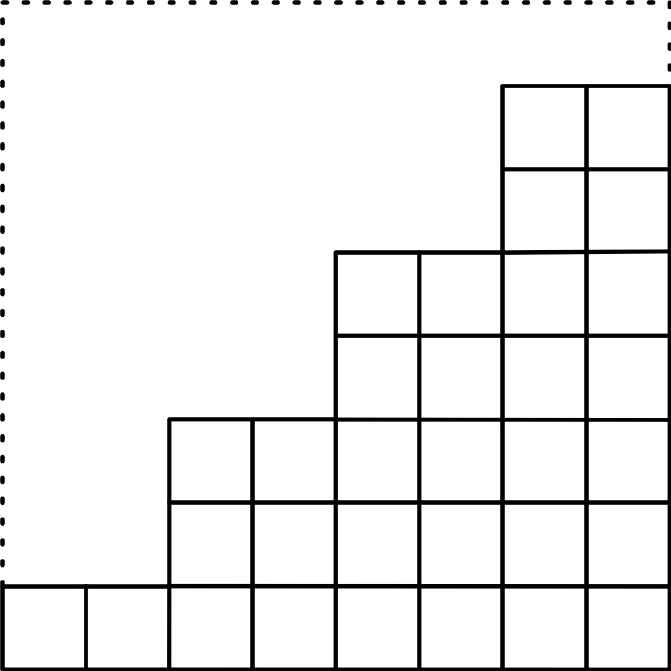}
\end{center}
\caption{Ciucu's way of halving the 8-by-8 square.}
\label{fig:ciucu}
\end{figure}

Henry Cohn~\cite{Co} proved that the function sending $n$ to $f(n)$
is uniformly continuous under the 2-adic metric
and thus extends to a function defined on all of $\Z$
and indeed all of $\Z_2$;
moreover, he showed that this extension satisfies
\begin{equation}
f(-1-n) = \left\{\begin{array}{rl}
 f(n) & \mbox{when $n$ is congruent to 0 or 3 (mod 4)}, \\
-f(n) & \mbox{when $n$ is congruent to 1 or 2 (mod 4)}.
\end{array} \right.
\end{equation}

Barkley and Liu~\cite{BL} have recently proved results
about 2-divisibility for the number of perfect matchings of a graph,
including as a special case the number of domino tilings of a rectangle,
but there is more refined work still to be done along the lines of Cohn's paper.
For instance, the mod 8 residue of 
the number of domino tilings of the $2n$-by-$(2n+2)$ rectangle
appears to depend only on the mod 4 residue of $n$;
the same goes for the number of domino tilings of the $2n$-by-$4n$ rectangle.

In this article we extend the discussion
to other sorts of tiles, specifically, tetrominos. 
A {\em tetromino} is a connected subset of the grid
that is a union of four grid-squares,
just as a domino is a union of two grid-squares.
Up to symmetry, there are five kinds of tetrominos: straight tetrominos, 
skew tetrominos, L-tetrominos, square tetrominos, and T-tetrominos.
They are shown in Figure~\ref{fig:the-six}, preceded by the domino.
These six tiles can be placed on a square grid in
2, 2, 4, 8, 1, and 4 translationally-inequivalent ways, respectively
(where rotations and reflections are permitted).
These are the sorts of tiles considered in this article.
({\em Trominos} -- unions of three grid-squares -- 
will be considered elsewhere.)

\begin{figure}[h]
\begin{center}
\includegraphics[width=4.5in]{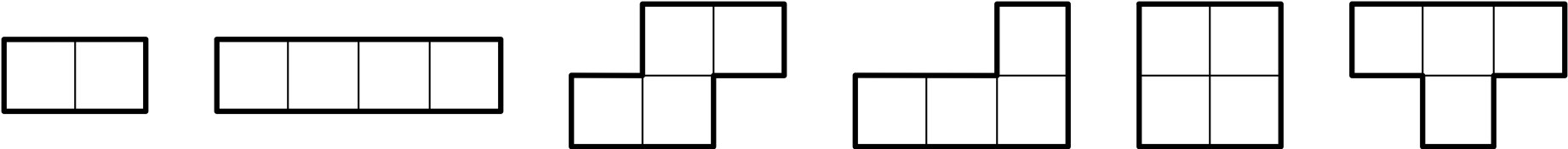}
\end{center}
\caption{A domino, a straight tetromino, a skew tetromino,
an L-tetromino, a square tetromino, and a T-tetromino.}
\label{fig:the-six}
\end{figure}

\end{section}

\begin{section}{Skew and straight tetrominos}
\label{sec:skewstraight}

\begin{figure}[h]
\begin{center}
\includegraphics[width=2in]{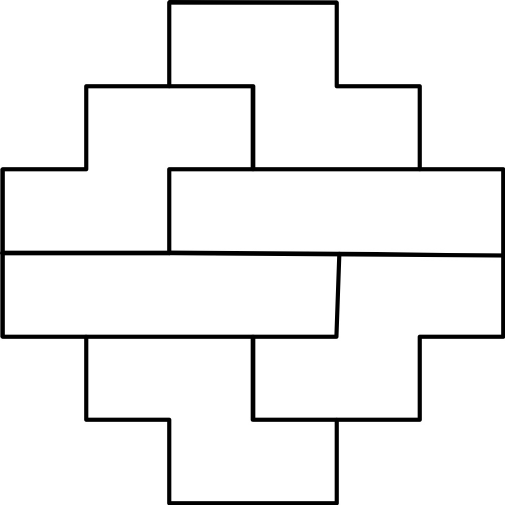}
\end{center}
\caption{Tiling the Aztec diamond of order 3 with skew and straight tetrominos.}
\label{fig:olympiad}
\end{figure}

I'll start with a warm-up puzzle that's 
roughly at the level of a math olympiad:
Prove that an Aztec diamond of order $n$ can be tiled 
by skew and straight tetrominos
(as shown in Figure~\ref{fig:olympiad} for $n=3$) 
only if $n$ is congruent to 0 or 3 (mod 4).

The puzzle can be solved using a valuation argument
(sometimes called a generalized coloring argument):
one can construct a mapping from the grid-cells to an appropriate abelian group
(a ``weight function'') and show that when $n$ is 1 or 2 (mod 4), 
the sum of the weights of the tiles
can't equal the sum of the weights of the region being tiled,
where the weight of a tile or a region being tiled
is the sum of the weights of the constituent cells.
Readers who are already familiar with this technique
might enjoy the challenge of attempting to solve the problem purely mentally.

\end{section}

\begin{section}{Dominos and square tetrominos}
\label{sec:domsquare}

In this section we use dominos and square tetrominos.
Thus an Aztec diamond of order 1 
(better known as the 2-by-2 square) can be tiled in 3 ways: 
with two horizontal dominos, two vertical dominos, or a single square tetromino.
The Aztec diamond of order 2 can be tiled in $2^{(2)(3)/2} = 8$ ways
using dominos, and can be tiled in an additional 11 ways
if one or more square tetrominos are included,
as shown in Figure~\ref{fig:eleven}.
Thus there are a total of $8+11=19$ ways
to tile an Aztec diamond of order 2
using dominos and square tetrominos.

\bigskip

\begin{figure}[h]
\begin{center}
\includegraphics[scale=0.4]{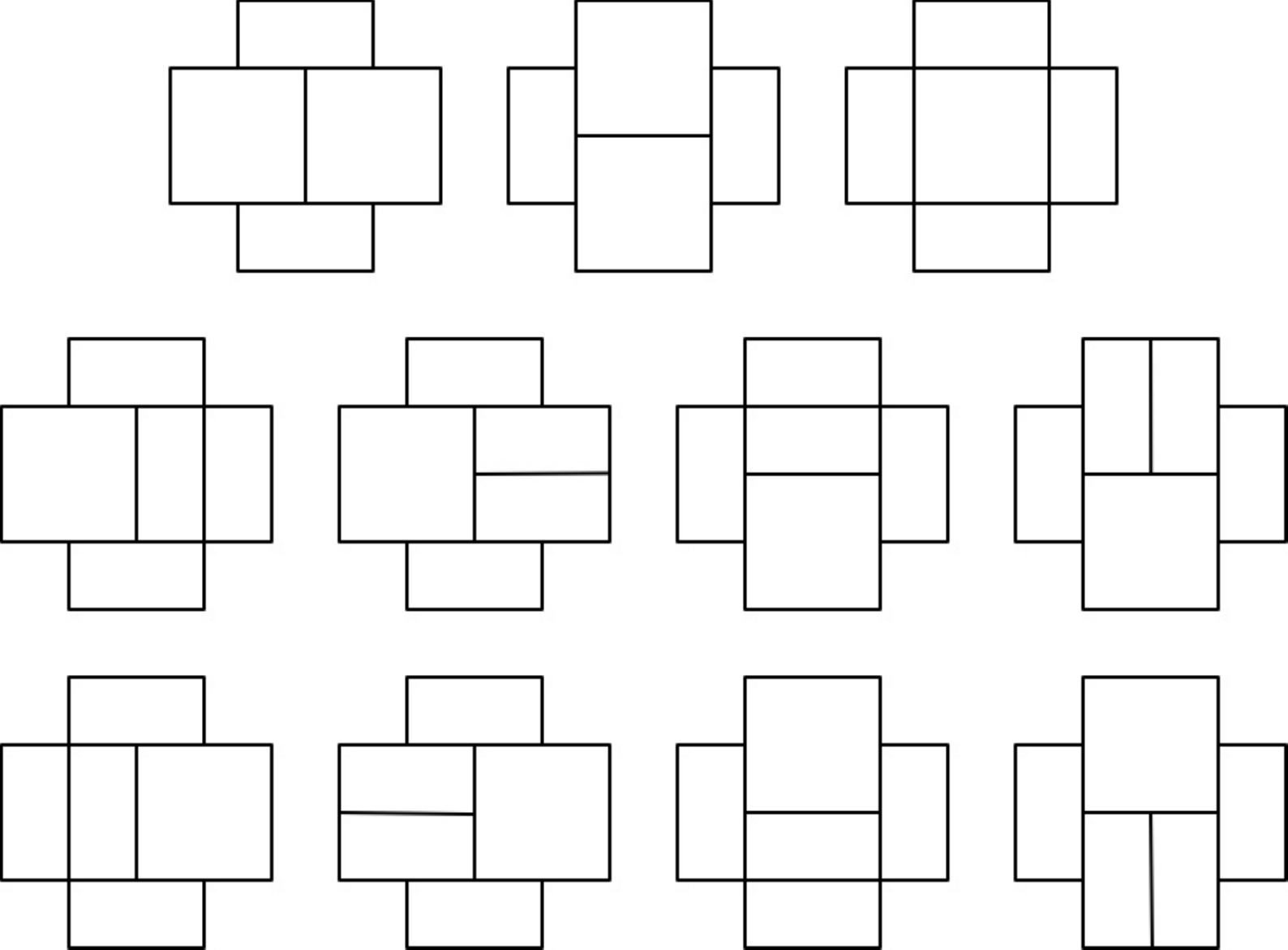}
\end{center}
\caption{Tiling the Aztec diamond of order 2
with dominos and at least one square tetromino.}
\label{fig:eleven}
\end{figure}

Define $M(n)$ (with $n \geq 0$) as 
the number of tilings of the Aztec diamond of order $n$
using dominos and square tetrominos. 
This is \href{https://oeis.org/A356512}{A356512}.
Trivially we have $M(0) = 1$ (since the Aztec diamond of order 0 is empty)
and we have already seen that $M(1) = 3$ and $M(2) = 19$. 
Figure~\ref{fig:Mtable} shows the terms of the sequence $M(n)$ 
for $n$ ranging from 0 to 12, 
computed using a program written by David desJardins.

The reader may wish to pause here to consider the problem of showing that $M(n)$ is always odd;
a solution will be given in section~\ref{sec:cong}.

\begin{figure}[h]
$$
\begin{array}{r|r}
n & M(n) \ \ \ \ \ \ \ \ \ \ \ \ \ \ \ \ \ \ \ \\
\hline
0 & 1 \\
1 & 3 \\
2 & 19 \\
3 & 293 \\
4 & 10917 \\
5 & 996599 \\
6 & 222222039 \\
7 & 121552500713 \\
8 & 162860556763865 \\
9 & 535527565429290907 \\
10 & 4318205059450240425083 \\
11 & 85475498697714319842817853 \\
12 & 4151186175463797888945512144221 
\end{array}
$$
\caption{Enumeration of tilings of Aztec diamonds
using dominos and square tetrominos.}
\label{fig:Mtable}
\end{figure}

These numbers grow quadratic-exponentially as a function of $n$,
and I have no conjectural formula for the $n$th term, 
nor a conjectural recurrence relation for the sequence, 
nor any efficient method of computing terms. 
Nonetheless, something systematic is going on. 
I have already mentioned that all the terms are odd. 
Taking this observation further, one notices that 
the numbers’ residues mod 4 are 
$$1, 3, 3, 1, 1, 3, 3, 1, 1, 3, 3, 1, 1,$$
the residues mod 8 are
$$1, 3, 3, 5, 5, 7, 7, 1, 1, 3, 3, 5, 5,$$
and the residues mod 16 are
$$1, 3, 3, 5, 5, 7, 7, 9, 9, 11, 11, 13, 13.$$

\bigskip

{\bf Conjecture 1}: For all $k \geq 1$, 
the mod $2^k$ residue of $M(n)$ is periodic with period dividing $2^k$. 
That is, $2^k$ divides $M(n+2^k) - M(n)$ for all $k,n$.

\bigskip

I tried to prove this conjecture by reducing it to an assertion 
about alternating-sign matrices but I was unsuccessful.

Note that if the conjecture is true
then $M(n) \equiv n + 1 + (1 + (-1)^{n + 1})/2$ (mod 8).
This congruence might also hold mod 16
but it certainly cannot hold mod $2^k$ for all $k$,
since that would require that $M(n)$ actually equals
$n + 1 + (1 + (-1)^{n + 1})/2$, which is clearly not the case for $n \geq 2$.
And indeed $M(2) = 19 \not\equiv 3$ (mod 32).

A deeper consequence of Conjecture 1 is that 
the function sending $n$ to $M(n)$ is 2-adically continuous.
Moreover, the function appears to satisfy a kind of symmetry
analogous to the functional equation (1) 
mentioned at the end of section~\ref{sec:intro}.

\bigskip

{\bf Conjecture 2}: For all $k \geq 1$,
if $n + n' \equiv -3$ (mod $2^k$) then $M(n) + M(n') \equiv 0$ (mod $2^k$).

\bigskip

That is, if one extends $M: \N \rightarrow \N$
to the 2-adic function $\widehat{M}: \Z_2 \rightarrow \Z_2$,
one has $\widehat{M}(-3-n) = - \widehat{M}(n)$.

Although in this article I am limiting myself 
to discussion of tilings of Aztec diamonds,
I have also looked at tilings of other regions 
using dominos and square tetrominos,
and the same phenomenon of 2-adic continuity arises fairly broadly there.
For instance, for the $2n$-by-$2n$ square, the $2n$-by-$(2n+2)$ rectangle, and the $2n$-by-$4n$ rectangle,
the number of tilings with dominos and square tetrominos
always seems to be congruent to $2n+1$ mod 8.

\end{section}

\begin{section}{Skew tetrominos and square tetrominos}
\label{sec:skewsquare}

In~\cite{Pr} I considered tilings of Aztec diamonds
by skew tetrominos and square tetrominos. 
If we require that all skew tetrominos be horizontal, 
interesting numerical patterns appear.
(Of course we would get the same result
if we required that all skew tetrominos be vertical.)
In this section we allow horizontal skew tetrominos and square tetrominos 
as seen in Figure~\ref{fig:tetra}, which depicts 
all six tilings of the Aztec diamond of order 3 
using square tetrominos and horizontal skew tetrominos.

\begin{figure}[h]
\begin{center}
\includegraphics[scale=0.35]{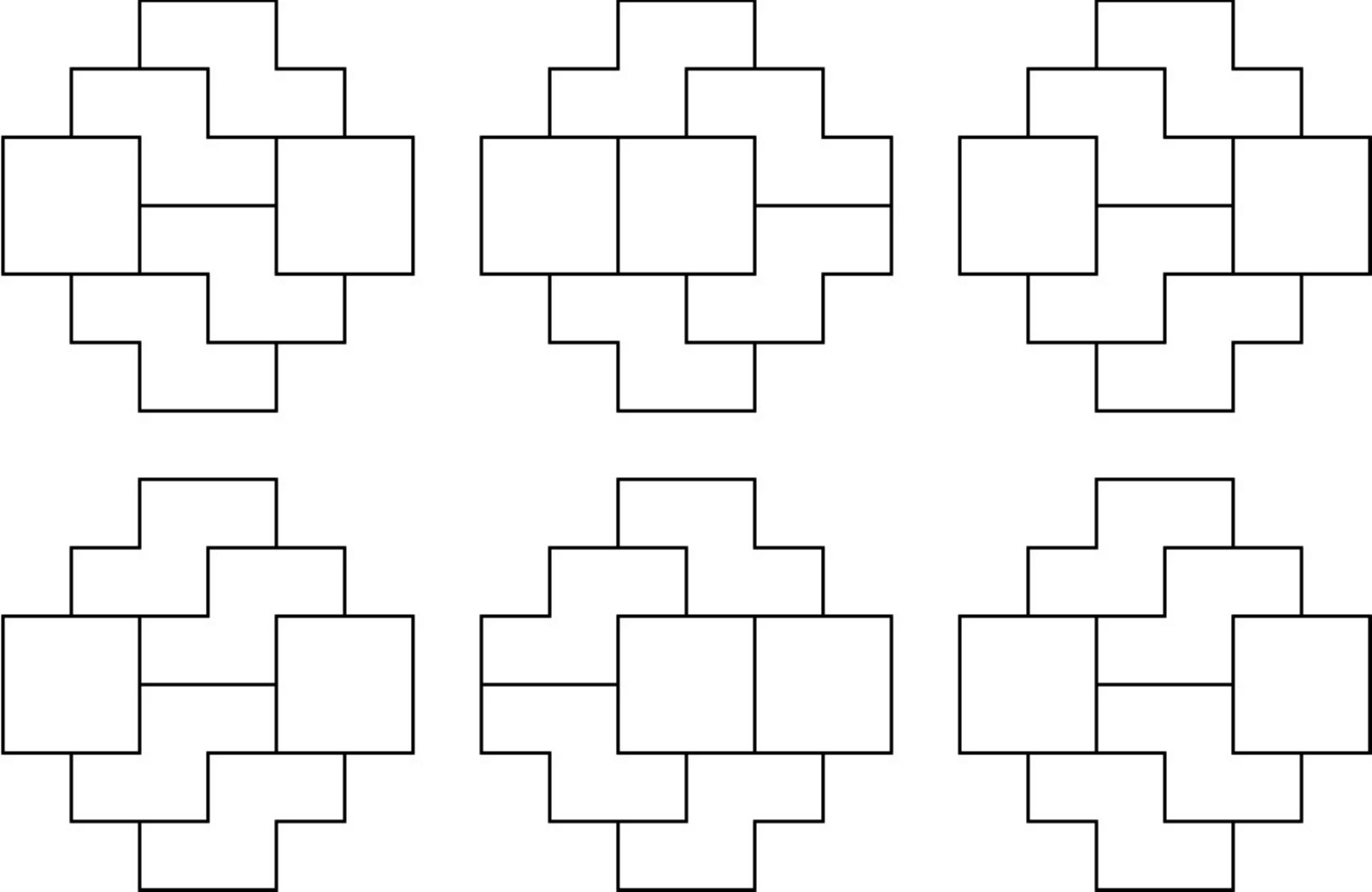}
\end{center}
\caption{Tiling the Aztec diamond of order 3
with horizontal skew tetrominos and square tetrominos.}
\label{fig:tetra}
\end{figure}

\begin{figure}[h!]
$$
\begin{array}{r|r}
n & L(n) \ \ \ \ \ \ \ \ \ \ \ \ \ \ \ \ \ \ \ \\
\hline
0 & 1 \\
1 & 1 \\
2 & 2 \\
3 & 6 \\
4 & 40 \\
5 & 364 \\
6 & 7904 \\
7 & 226152 \\
8 & 15835008 \\
9 & 1439900880 \\
10 & 324189571584 \\
11 & 94080051207136 \\
12 & 68041472016287744 \\
13 & 63145927127133361600 \\
14 & 146637148542938673930240 \\
15 & 435697213021432661980535936
\end{array}
$$
\caption{Enumeration of tilings of Aztec diamonds
using horizontal skew tetrominos and square tetrominos.}
\label{fig:Ltable}
\end{figure}

Define $L(n)$ (with $n \geq 0$) as 
the number of tilings of the Aztec diamond of order $n$
using horizontal skew tetrominos and square tetrominos.
This is \href{https://oeis.org/A356513}{A356513}.
Trivially we have $L(0) = 1$ and $L(1) = 1$. 
Figure~\ref{fig:Ltable} shows the terms 
of the sequence $L(n)$ for $n$ ranging from 0 to 15, 
again computed using the program written by David desJardins.

The sequence grows quadratic-exponentially,
and once again, I have no conjectural formula,
but as before there are patterns that call out for explanation. 
Noticing that all but the first two terms are even, 
one might naturally think to look at the multiplicity of 2 
in the prime factorization of $L(n)$, obtaining the sequence
0,0,1,1,3,2,5,3,7,4,9,5,11,6,13,7,\dots,
which (once we throw out the initial 0)
we recognize as an interspersal of the arithmetic progressions 
0,1,2,3,4,5,6,7,\dots  and 1,3,5,7,9,11,13,\dots.

\bigskip

{\bf Conjecture 3}: 
For $n \geq 1$,
the multiplicity of 2 in the prime factorization of $L(n)$ 
is $n-1$ if $n$ is even and $(n-1)/2$ if $n$ is odd.

\bigskip

Going further, 
let $L_0(m) = L(2m) / 2^{2m-1}$ and $L_1(m) = L(2m-1) / 2^{m-1}$, 
so that (if Conjecture 3 holds) 
$L_0(m)$ and $L_1(m)$ are odd integers for all $m$. These two new sequences
are shown in Figure~\ref{fig:JKtable}.

\begin{figure}[h!]
$$
\begin{array}{r|r}
m & L_0(m) \ \ \ \ \ \ \ \ \ \ \ \\
\hline
1 & 1 \\
2 & 5 \\
3 & 247 \\
4 & 123711 \\
5 & 633182757 \\
6 & 33223375007953 \\
7 & 17900042546745443595 \\
\ & \ {}
\end{array}
\ \ \ \begin{array}{r|r}
m & L_1(m) \ \ \ \ \ \ \ \ \ \ \ \ \ \ \ \\
\hline
1 & 1 \\
2 & 3 \\
3 & 91 \\
4 & 28269 \\
5 & 89993805 \\
6 & 2940001600223 \\
7 & 986655111361458775 \\
8 & 3403884476729942671722937
\end{array}
$$
\caption{Values of $L_0(m)$ and $L_1(m)$.}
\label{fig:JKtable}
\end{figure}

The mod 4 residues of the $L_0$ sequence go $1, 1, 3, 3, 1, 1, 3$
while those of the $L_1$ sequence go $1, 3, 3, 1, 1, 3, 3, 1$.
That's not much evidence to go on,
so perhaps it would be prudent not to make a conjecture,
but I choose to be hopeful.

\bigskip

{\bf Conjecture 4}: For all $k \geq 0$, the mod $2^k$ residue of $L_0(m)$ 
is periodic with period dividing $2^k$. Likewise for $L_1(m)$.

\bigskip

We do not observe such patterns in the numbers of tilings
when both horizontal and vertical skew tetrominos are allowed
along with square tetrominos as in~\cite{Pr}.
More specifically, if we count tilings of Aztec diamonds
in which we are permitted to use all four kinds of skew tetrominos
as well as square tetrominos, the resulting sequence, taken mod 4, goes
1, 0, 0, 0, 0, 0, 2, 2, 0, 0, 0, 0, \dots;
if there is a period here, and it is a power of 2,
it must be at least 16.

The prime $p=2$ appears to be special
for the enumerative problems I described above;
looking at the $M$ and $L$ sequences 
mod 3 or mod 5 yields no discernible patterns.

\end{section}

\begin{section}{Assorted congruential problems}
\label{sec:cong}

For each of the sixty-three nonempty subsets of 
the set of six tiles shown in Figure~\ref{fig:the-six},
we can ask in how many ways 
it is possible to tile the Aztec diamond of order $n$
using only tiles from that set,
allowing translations, rotations, and reflections of tiles.
These are the enumerative problems considered in this section.

(One could expand the set of tiling problems
by distinguishing between different orientations of the tiles,
as was done in the preceding section
where we permitted horizontal skew tetrominos
but forbade vertical skew tetrominos;
since there are $2+2+4+8+1+4=21$ different tiles up to translation,
we would obtain over two million different problems,
and even if we mod out the $2^{21}-1$ problems by a dihedral action of order 8,
that is still too many problems to consider exhaustively.
One that appears to be interesting is discussed at the end of this section.)

In each of the sixty-three cases, I used the aforementioned program
to count the tilings of the Aztec diamond of order $n$,
with $n$ going from 1 to 8, using the allowed tiles.
Although no 2-adic continuity phenomena arose from these experiments,
there were definite patterns in the parity,
and in a few cases there were congruence patterns modulo higher powers of 2.
Here I will adopt a six-bit code to represent
the sixty-three tiling problems,
in which the six successive bits (from left to right) equal 1 or 0
according to whether or not dominos, straight tetrominos, skew tetrominos, 
L-tetrominos, square tetrominos, and T-tetrominos are allowed.
For instance, the case treated in section~\ref{sec:skewstraight},
in which only straight tetrominos and skew tetrominos are allowed
(see Figure~\ref{fig:olympiad}), would be assigned the code 011000;
the case treated in section~\ref{sec:domsquare},
in which only dominos and square tetrominos are allowed
(see Figure~\ref{fig:eleven}), would be assigned the code 100010;
and the case of unconstrained skew and square tetrominos
(briefly discussed in section~\ref{sec:skewsquare}) 
would be assigned the code 001010.

In one-third of the 63 cases,
I observed that for all $n$ between 1 and 8,
the number of tilings of the Aztec diamond of order $n$ is even.
These were the cases associated with the six-bit codes
001001, 001100, 001101, 011001, 011100, 011101, 100001,
100100, 100101, 101000, 101001, 101100, 101101, 110000,
110001, 110100, 110101, 111000, 111001, 111100, and 111101.

Presumably some (perhaps all) of these examples can be resolved 
by showing that there are no tilings that are invariant
under the full dihedral group, since in that case all orbits 
would contain an even number of tilings.

Three of the 21 cases were especially interesting.
In case 011100, all terms were divisible by 8;
in case 100001, all terms after the first were divisible by 8;
and in case 110001, all terms were congruent to 2 (mod 4).

There were also four cases in which I observed that
the number of tilings of the Aztec diamond of order $n$ is even
for all $n$ between 2 and 8
(with the number of tilings being the odd number 1 in the case $n=1$).
These were the cases associated with the six-bit codes
001010, 001110, 011010, and 011110.

In the cases 001101, 100001, 100011, and 111000
it appears that the exponent of 2 in the number of tilings
may be going to infinity with $n$,
though with such scant evidence it would be rash
to place too much faith in this guess. 

Additionally, there is one case in which the number of tilings 
of the Aztec diamond of order $n$ is always odd,
namely, tilings using only dominos and square tetrominos.
Indeed, if we assign each tiling weight $(-1)^s$ where $s$ is the number of square tetrominos, 
I claim that the sum of the weights is 1. 
We can prove this using a sign-reversing involution 
that scans through the tiling in some fashion in search of a 2-by-2 block 
that is tiled either with a square tetromino or with two vertical dominos 
and switches between the two possibilities. 
The fixed points of this involution are tilings that use only horizontal dominos, 
and there is just one of those.
[NOTE: A referee pointed out that the preceding proof is incorrect
and suggested a way to fix it. This change was implemented in the
final, published version of the article at
\url{https://math.colgate.edu/~integers/x30/x30.pdf}.]

Finally, leaving the small world of the $2^6-1$ problems
and dipping our toe into the big world of the $2^{21}-1$ problems,
we consider tilings of the Aztec diamond of order $n$
by dominos and horizontal straight tetrominos.
This is \href{https://oeis.org/A356523}{A356523}, and begins
1, 2, 11, 209, 12748, 2432209, 1473519065, \dots.
It appears that the number of tilings is 
even when $n \equiv 1$ (mod 3) and odd otherwise;
this has been verified for $1 \leq n \leq 16$.

\end{section}

\begin{section}{Reduction to perfect matchings}
\label{sec:matchings}

The $L$ sequence from section~\ref{sec:skewsquare}
has an interpretation in terms of perfect matchings.
To see why, suppose we have a tiling of the Aztec diamond of order $n$
using horizontal skew tetrominos and square tetrominos.
Dividing each tetromino into two horizontal dominos
gives us a tiling of the Aztec diamond by horizontal dominos,
but it is easy to see that there is exactly one such tiling (call it $T$).
Hence each tetromino is obtained by gluing together two dominos in $T$.
That is, the tetromino tilings correspond to perfect matchings in the graph
whose vertices correspond to the dominos in $T$
with an edge joining two vertices if the corresponding dominos
form a horizontal skew tetromino or square tetromino.
It is not hard to see that this graph is similar
to the $n$-by-$n$ square except that the diagonal has been ``doubled'';
for instance, the right panel of Figure~\ref{fig:double-double}
shows the graph for $n=4$.

\begin{figure}[h]
\begin{center}
\includegraphics[width=4in]{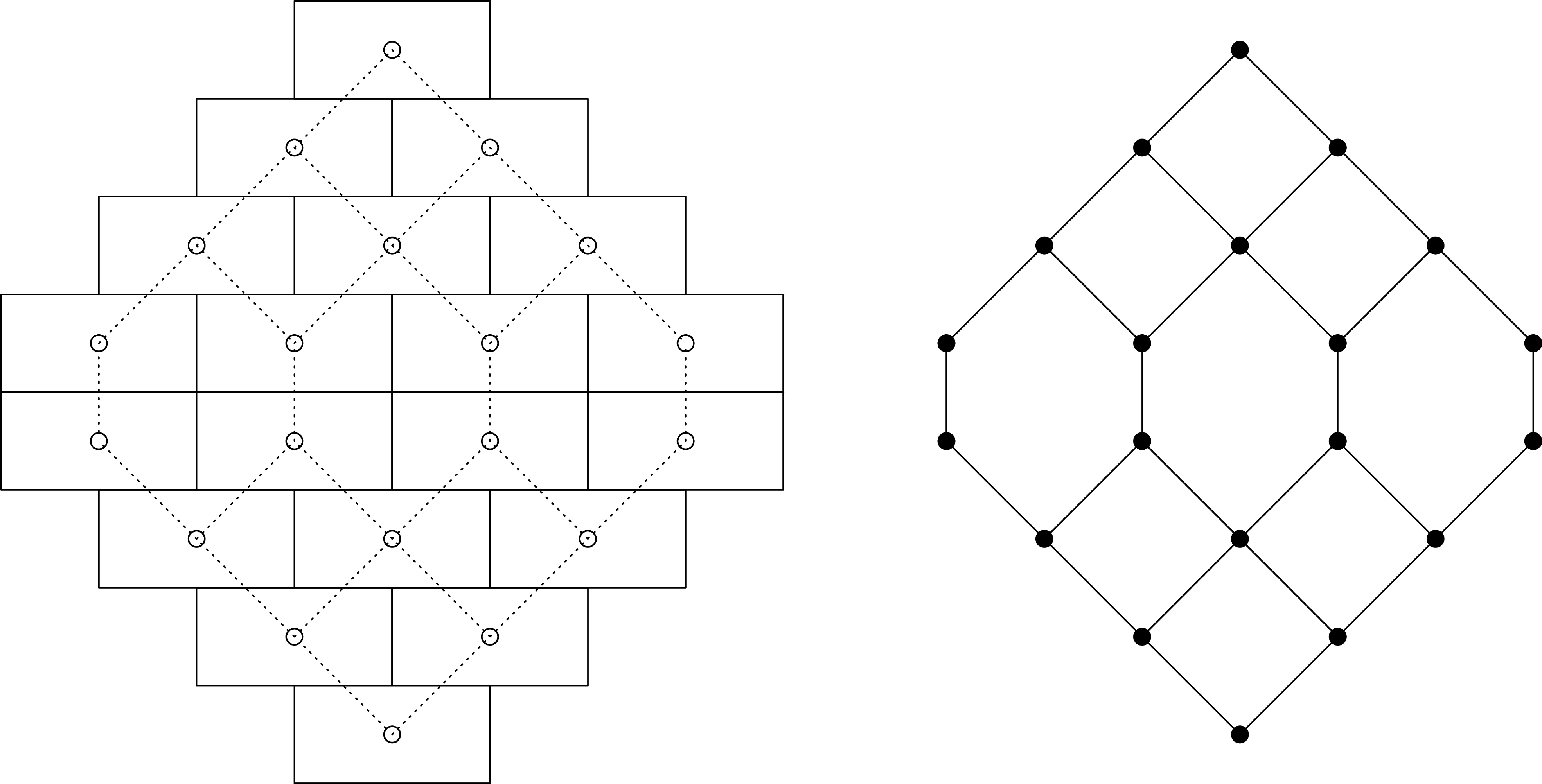}
\end{center}
\caption{Deriving the square graph with doubled diagonal.}
\label{fig:double-double}
\end{figure}

A similar analysis can be applied to tilings of Aztec diamonds
using horizontal skew tetrominos and horizontal straight tetrominos.
In this case the Aztec diamond splits 
into two non-interacting halves (top half and bottom half),
each of which can be tiled independently of the other, 
and the tilings of either half correspond 
to perfect matchings of a triangle graph
as shown in Figure~\ref{fig:triangles}.
Thus the number of such tetromino tilings of the Aztec diamond of order $n$
is equal to the square of the $n$th term of sequence 
\href{https://oeis.org/A071093}{A071093}.
Studying the first 25 terms, I find that
the sequence seems to have 2-adic properties of its own.
The largest power of 2 dividing the $n$th term of the sequence A071093
appears to be $\lfloor n/2 \rfloor$,
and the 2-free part appears to satisfies 2-adic continuity:
for instance, its value mod 16 seems to be determined by $n$ mod 16.

\begin{figure}[h]
\begin{center}
\includegraphics[width=4in]{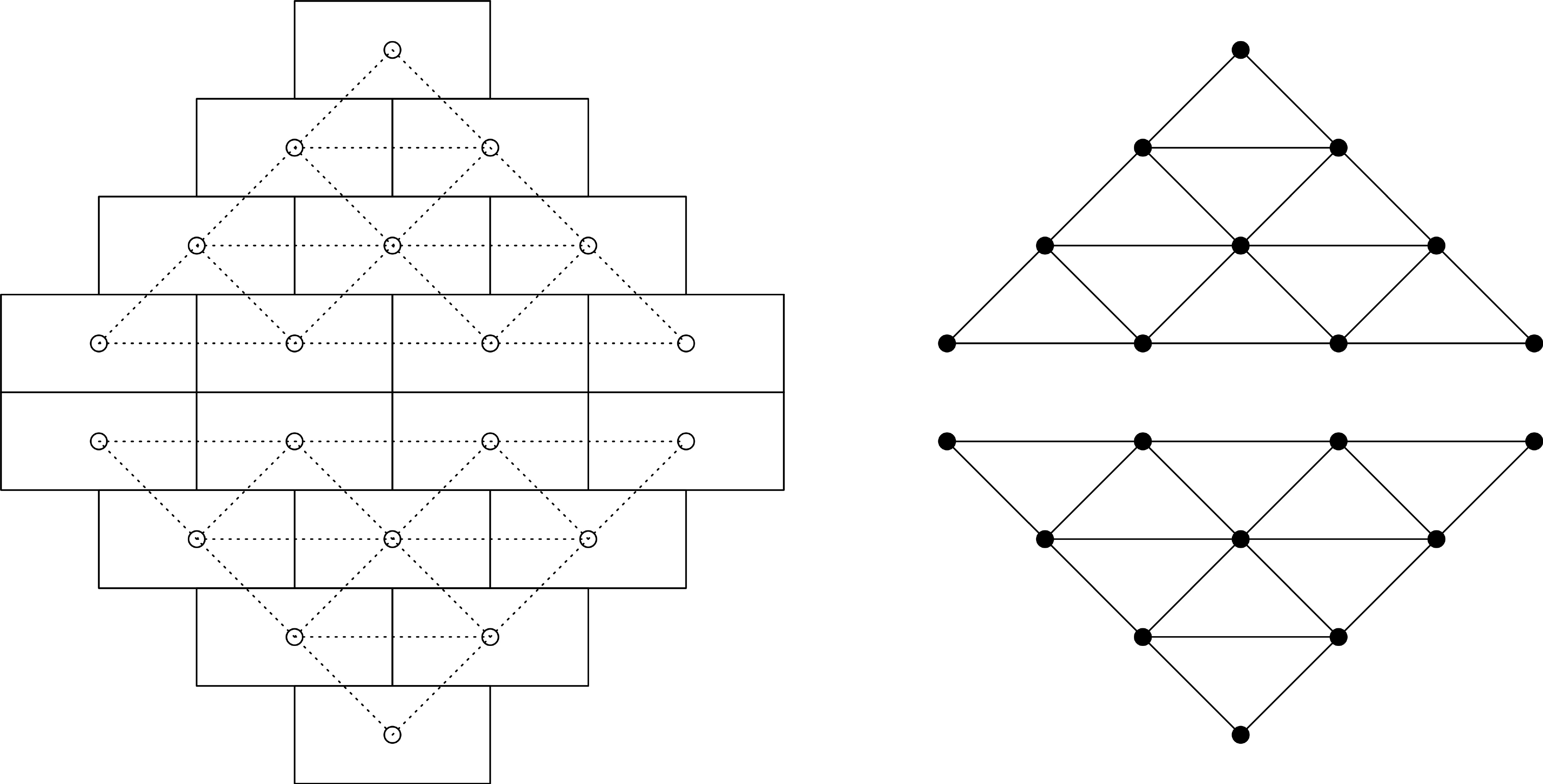}
\end{center}
\caption{Deriving the double triangle graph.}
\label{fig:triangles}
\end{figure}

What if we superimpose the two graphs, 
obtaining the graph shown at the right half of
Figure~\ref{fig:superimpose}?
This is equivalent to tiling an Aztec diamond using
horizontal skew tetrominos, horizontal straight tetrominos, and square tetrominos.
Then, counting the tilings, we obtain the integer sequence
1, 2, 10, 116, 3212, 209152, 32133552, 11631456480, 9922509270288, 
19946786274879008, 94492874103638971552, 1054865198752147761744448,
\dots.
This is \href{https://oeis.org/A356514}{A356514}.
It appears that the number of tilings is divisible
by $2^{\lfloor n/2 \rfloor}$.

\begin{figure}[h]
\begin{center}
\includegraphics[width=4in]{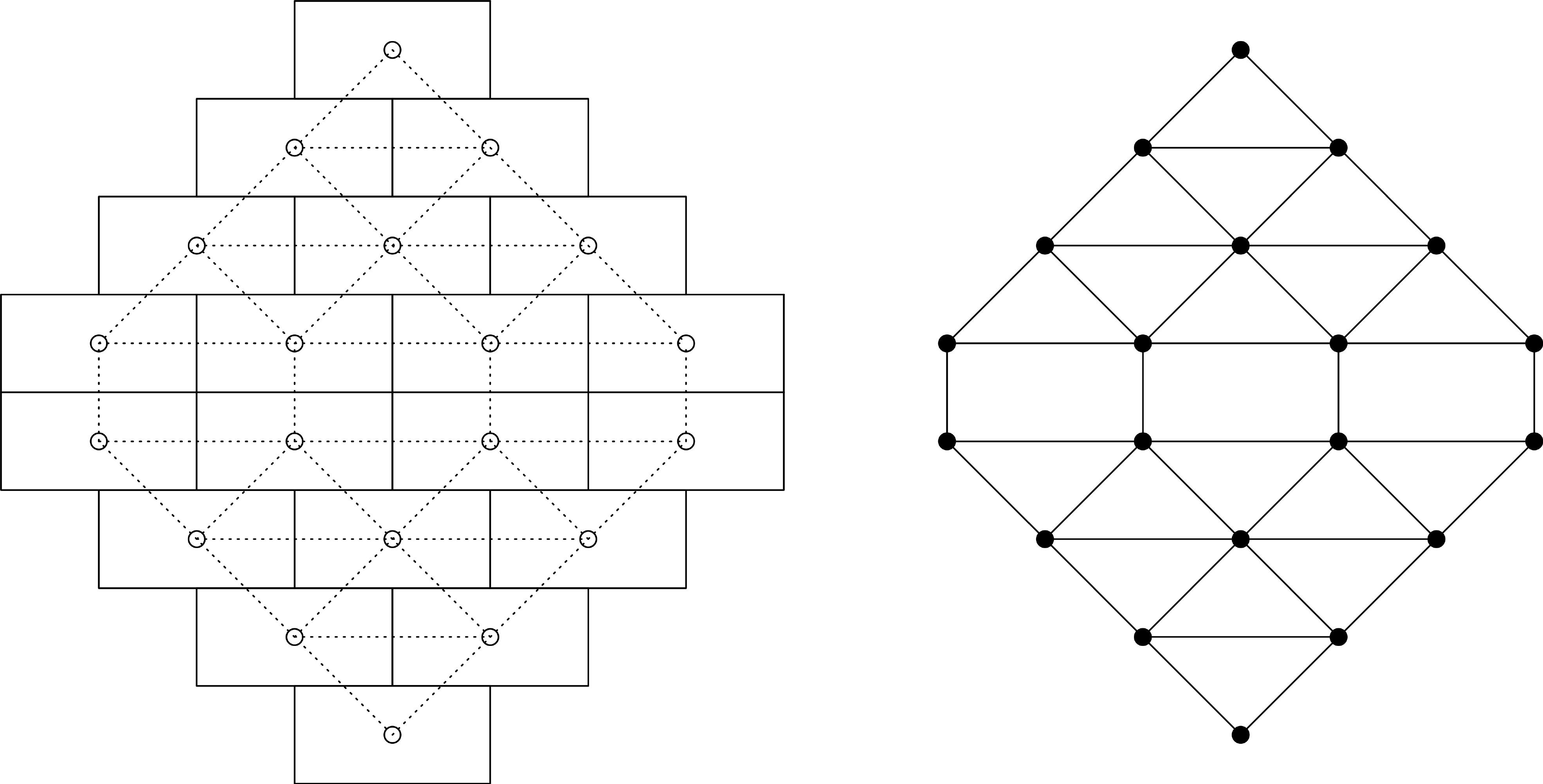}
\end{center}
\caption{Superimposing Figures~\ref{fig:double-double} and~\ref{fig:triangles}.}
\label{fig:superimpose}
\end{figure}

\end{section}

\begin{section}{Some thoughts}
\label{sec:thoughts}

The articles of Lovasz~\cite{Lo}, Ciucu~\cite{Ci}, 
Pachter~\cite{Pa}, and Barkley and Liu~\cite{BL}
give ways to find the largest power of 2
that divides the number of perfect matchings of a graph.
This should provide traction for Conjecture 3,
since we saw in section~\ref{sec:matchings} that
the $L$ sequence has an interpretation in terms of perfect matchings
of certain graphs.
Graphs of this kind appear in the paper of Ciucu~\cite{Ci};
in particular, his Lemma 1.1 shows that
the number of perfect matchings is divisible by $2^{\lfloor n/2 \rfloor}$.
By bringing ideas from Pachter~\cite{Pa}, one might be able to prove Conjecture 3,
as well as some of the other 2-divisibility conjectures from this article.

The only work I know of that provides detailed 2-adic information 
about the 2-free part of numbers that count tilings
is the work of Cohn~\cite{Co}.
Cohn's approach presupposes the existence of an exact formula
(in Cohn's case, an explicit product of algebraic integers);
perhaps something similar can be done for perfect matchings
of the square graph with doubled diagonal,
yielding a proof of Conjecture 4.

Conjectures 1 and 2 seem harder.
The product formula exploited by Cohn 
was discovered by Temperley and Fischer~\cite{TF} 
and independently by Kasteleyn~\cite{Ka} at about the same time;
those researchers made use of the fact that,
just as determinants and Pfaffians of matrices can be expressed as sums of terms
associated with perfect matchings of the set of rows and columns,
one can conversely express the number of perfect matchings of a planar graph
in terms of the determinant or Pfaffian of an associated matrix.
I know of no way of recast the $M$ sequence 
as enumerating perfect matchings of graphs.
However, it is easy to recast the $M$ sequence as enumerating
perfect matchings of certain hypergraphs.
Can any of the existing notions of hyperdeterminants be brought to bear?
Perhaps a reading of~\cite{GKZ} would suggest possible approaches.

Kuperberg's elegant solution~\cite{Ku} 
to the alternating sign matrix conjecture
exploits the power of the Yang-Baxter equation in statistical mechanics.
It's possible that tools for analyzing the new problems described in this article
will be found in the existing literature at 
the interface between algebra and statistical mechanics.

In any case, inasmuch as Conjectures 1 and 2 are reminiscent of Cohn's work,
and inasmuch as Cohn's argument hinges on an exact product formula, 
one might hope that an exact formula of some kind
can be found for the $M$ sequence.  Such an exact formula would have other uses.
In~\cite{CEP} and~\cite{CLP}, Henry Cohn, Noam Elkies, Michael Larsen and I 
used exact enumeration results to prove 
concentration theorems for random tilings.
One might hope that the curious 2-adic phenomena discussed in this article
hint at the existence of algebraic machinery that
could be applied to the task of showing us
what random tilings associated with Conjecture 1
look like in the limit as size goes to infinity.
Preliminary experiments suggest that 
there is a ``frozen region'' near the boundary,
but I have no idea how far into the interior it extends.

\end{section}

\bigskip
\bigskip

\noindent
I thank David desJardins for the software that made this research possible,
and Noam Elkies for helpful comments.
Above all I thank Michael Larsen for his many years 
of friendship and mathematical camaraderie.

\bibliographystyle{abbrvnat}
\bibliography{sample}
\label{sec:biblio}

\end{document}